\documentclass{amsart}
\usepackage{geometry} 
\usepackage{url}
\usepackage{amsthm}
\usepackage{amssymb}
\usepackage{latexsym}
\usepackage{amsfonts}
\usepackage{amsmath}
\usepackage{amstext}
\usepackage{accents}
\usepackage{cancel}
\usepackage{eucal}
\usepackage{amscd}
\usepackage{hyperref}
\usepackage[shortalphabetic]{amsrefs}

\begin{document}

\title[Evolution of closed curves by length-constrained curve diffusion]{Length-constrained curve diffusion} 
\author[J. McCoy]{James McCoy*} \thanks{* Corresponding author}
\address{School of Mathematical and Physical Sciences, University of Newcastle}
\email{James.McCoy@newcastle.edu.au}
\author[G. Wheeler]{Glen Wheeler}
\address{Institute for Mathematics and its Applications, University of Wollongong}
\email{glenw@uow.edu.au}
\author[Y. Wu]{Yuhan Wu}
\address{Institute for Mathematics and its Applications, University of Wollongong}
\email{yw120@uowmail.edu.au}

\thanks{The research of the first and second authors was supported by Discovery Project grant DP150100375 of the Australian Research Council.  Some of this research was conducted while the first author was visiting the University of Wollongong.  The research of the third author was supported by a University of Wollongong Faculty of Engineering and Information Sciences Postgraduate research scholarship.}
\begin{abstract}
We show that any initial closed curve suitably close to a circle flows under length-constrained curve diffusion to a round circle in infinite time with exponential convergence.  We provide an estimate on the total length of time for which such curves are not strictly convex.  We further show that there are no closed translating solutions to the flow and that the only closed rotators are circles.
\end{abstract}

\maketitle

\section{Introduction} \label{S:intro}
\newtheorem{main}{Theorem}[section]
\newtheorem{pos}[main]{Proposition}
\newtheorem{embed}[main]{Proposition}
\newtheorem{transl8}[main]{Proposition}
\newtheorem{rot8}[main]{Proposition}
In a recent article \cite{W13} the second author considered the curve diffusion flow of closed plane curves.  This flow has the fundamental property that the signed enclosed area is preserved under the flow while the length of the evolving curve does not increase.  As such the flow provides a natural approach to the isoperimetric problem.  It is also natural therefore to consider a `dual' fourth order flow that preserves length of the evolving curve while the signed enclosed area does not decrease.  Such a flow may be obtained by including an appropriate globally-defined function of time in the flow speed.  Specifically, suppose $\gamma : \mathbb{S}^1 \times \left[ 0, T\right) \rightarrow \mathbb{R}^2$ evolves by the fourth-order curvature flow
\begin{equation} \label{E:theflow}
  \partial_t^\perp \gamma \left( x, t\right) = h\left( t\right) - k_{ss} \left( x, t\right)\mbox{,}
\end{equation}
where $k_{ss}$ denotes the second derivative of (scalar) curvature of $\gamma$ with respect to arc length $s$.  As usual for geometric flow problems, we need only specify the normal component $\partial_t^\perp \gamma$ of the evolution of $\gamma$; any tangential component corresponds to reparametrisations of the evolving curve (to obtain short-time existence of a family of solution curves we would fix a parametrisation by adding a specific tangential term to \eqref{E:theflow}).  To preserve length of the evolving curve $\gamma_t := \gamma \left( \cdot, t\right)$ we take
\begin{equation} \label{E:h}
  h\left( t\right) = - \frac{\int k_s^2 ds}{2\pi \omega} \mbox{,}
\end{equation}
where $\omega$ denotes the winding number of $\gamma_t$.  One may then obtain short-time existence of a solution to \eqref{E:theflow} using a standard fixed point argument between appropriate function spaces with suitably smooth functions $f(t)$ in place of $h\left( t\right)$.  For details of a similar fixed point argument (but for a second order flow) we refer the reader to \cite{M05}; for a discussion of the different approaches available to short-time existence for the regular curve diffusion flow see \cite{W13} and the references contained therein (adding $f\left( t\right)$ term causes no difficulty as it does not change the symbol of the differential operator).

Writing $I\left[ \gamma\right] = \frac{L^2}{4\pi A}$ for the isoperimetric ratio of $\gamma$ and $K_{\mbox{osc}}\left[ \gamma\right] = L \int \left( k - \overline{k}\right)^2 ds$ for the normalised oscillation of curvature, where $\overline{k} = \frac{\int k ds}{L}$ is the average curvature at time $t$, our main theorem is as follows

\begin{main} \label{T:main}
Suppose $\gamma_0: \mathbb{S}^1 \rightarrow \mathbb{R}^2$ is a regular smooth immersed closed curve with $A\left[ \gamma_0\right] >0$, length $L_0$ and 
$$\int k \, ds = 2 \pi \mbox{.}$$
There exists a constant $K^*>0$ such that if
$$K_{\mbox{osc}}\left[ \gamma_0\right] < K^* \mbox{ and } I\left[ \gamma_0\right] < \frac{4\pi^2}{4 \pi^2 - K^*}$$
then the length-constrained curve diffusion flow \eqref{E:theflow} with initial data $\gamma_0$ exists for all time and converges exponentially to a round circle with radius $\frac{L_0}{2\pi}$.
\end{main}

There is a couple of other related results that may be proved similarly as in \cite{W13}.  The first is an upper bound of the size of the set of times for which the curvature of a solution of \eqref{E:theflow} is not strictly positive. 

\begin{pos} \label{T:pos}
Suppose $\gamma: \mathbb{S}^1 \times \left[ 0, T\right) \rightarrow \mathbb{R}^2$ solves \eqref{E:theflow} and the assumptions of Theorem \ref{T:main}.  Then
$$\mathcal{L}\left\{ t \in \left[ 0, \infty \right) : k\left( \cdot, t\right) \not > 0 \right\} \leq \frac{L_0^2}{4\, \pi^3}\left( \frac{L_0^2}{4\, \pi} - A_0 \right)$$
where $A_0$ denotes the initial enclosed area and $k\left( \cdot, t\right) \not > 0$ means there exists an $x\in \mathbb{S}^1$ such that $k\left( x, t\right) \leq 0$.
\end{pos}

The estimate of Proposition \ref{T:pos} is optimal in the sense that for a circle (whose image is static under \eqref{E:theflow}) the right hand side is equal to zero.\\

The result that solutions with sufficiently small oscillation of curvature remain embedded for all time applies exactly as in \cite{W13}:

\begin{embed}
Any solution of \eqref{E:theflow} with initial embedded curve $\gamma_0$ satisfying the assumptions of Theorem \ref{T:main} remains embedded for all time.
\end{embed}

Our results on self-similar closed curves evolving under \eqref{E:theflow} are as follows:

\begin{transl8} \label{T:transl8}
Let $\gamma: \mathbb{S}^1\rightarrow \mathbb{R}^2$ be a smooth, closed, translating solution of \eqref{E:theflow} with circumference length $L_0$.  Then $\gamma\left( \mathbb{S}^1\right)$ is the stationary round circle of circumference $L_0$.
\end{transl8}

\begin{rot8} \label{T:rot8}
 Let $\gamma: \mathbb{S}^1\rightarrow \mathbb{R}^2$ be a smooth, closed, rotating solution of \eqref{E:theflow} with circumference length $L_0$.  Then $\gamma\left( \mathbb{S}^1\right)$ is a standard round circle of circumference $L_0$.
 \end{rot8}

Higher order geometric evolution problems have received increasing attention in
the last few years. Particular geometric fourth order equations occur in
physical problems and enjoy some interesting applications in mathematics.
We mention in particular for curves the curve diffusion flow and $L^2$-gradient flow of the
elastic energy, and for surfaces the surface diffusion and Willmore flows.  Fourth order flows with constraints have been considered for example in \cites{MWW11, MW}.  Relevant work on higher order flows of closed curves without boundary includes \cites{DKS02, EGBMWW14, GI99, PW16, W13, W15}.

The remainder of this article is organised as follows.  In Section \ref{S:setup} we define notation, state some key tools to be used in our analysis and provide a bound on $\left| h\left( t\right)\right|$ via interpolation and the evolution equations for the various geometric quantities we will need.  In Section \ref{S:estimates} we focus on estimating the oscillation of curvature.  With these results in hand in Section \ref{S:global} we prove long-time existence and convergence to circles under the stated conditions, completing the proof of Theorem \ref{T:main} and \ref{T:pos}.  Finally in Section \ref{S:selfsim} we prove the results on self-similar solutions.

\section{Preliminaries} \label{S:setup}
\newtheorem{PSW}{Lemma}[section]
\newtheorem{interp}[PSW]{Lemma}
\newtheorem{interp2}[PSW]{Lemma}
\newtheorem{hbound}[PSW]{Corollary}
\newtheorem{evlneqns}[PSW]{Lemma}
\newtheorem{ste}[PSW]{Theorem}
Let $\gamma_0 : \mathbb{R} \rightarrow \mathbb{R}^2$ be a (suitably)
smooth embedded (or immersed) regular curve.  We say $\gamma$ is \emph{periodic} with period $P$ if there exists a vector $V\in \mathbb{R}^2$ and a positive number $P$ such that, for all $m\in \mathbb{N}$,
$$\gamma\left( x + P \right) = \gamma\left( x\right) + V \mbox{ and } \partial_x^m \gamma \left( x + P\right) = \partial_x^m \gamma\left( x\right) \mbox{.}$$
Here $\partial_x^m$ denotes the $m$th iterated derivative of $\gamma$.  If $V=0$ then $\gamma$ is closed and we may rewrite $\gamma: \mathbb{S}^1 \rightarrow \mathbb{R}^2$.  The length of $\gamma$ is
$$L\left[ \gamma\right] = \int_0^P \left| \gamma'\left( u \right) \right| du$$
and the signed enclosed area is 
$$A\left[ \gamma \right] = - \frac{1}{2}\int_0^P \left< \gamma, \nu\right> \left| \gamma' \right| du \mbox{,}$$
where $\nu$ is a unit normal vector field on $\gamma$.  Throughout this article we will keep our evolving curves $\gamma$ parametrised by arc length $s$.\\

We will frequently use the following Poincar\'{e}-Sobolev-Wirtinger `[PSW]' inequalities.  For proofs of these see for example Appendix A of \cite{PW16}.

\begin{PSW} \label{T:PSW}
Suppose $f: \mathbb{R}\rightarrow \mathbb{R}$ is absolutely continuous and periodic with period $P$.  Then if $\int_0^P f\left( x\right) dx =0$ we have
\begin{enumerate}
  \item[\textnormal{(i)}]
$$\int_0^P f^2\left( x\right) dx \leq \frac{P^2}{4\pi^2} \int_0^P \left| f_x\left( x\right) \right|^2 dx$$
with equality if and only if $f\left( x\right) = a\, \sin \left( \frac{2\, \pi\, x}{P} + b \right)$ for arbitrary constants $a$ and $b$;
  \item[\textnormal{(ii)}]
$$\left\| f \right\|_\infty^2 \leq \frac{P}{2\pi} \int_0^P \left| f_x\left( x\right) \right|^2 dx \mbox{.}$$
\end{enumerate}
\end{PSW}

We also need some interpolation inequalities from \cite{DKS02}.   We first need to set up some notation.  For normal tensor fields $S$ and $T$ we denote by $S \star T$
 any linear combination of $S$ and $T$.  In our setting, $S$ and $T$ will be simply curvature $k$ or its arc length derivatives.  Denote by $P_n^m\left( k\right)$ any linear combination of terms of type $\partial_s^{i_1} k \star \partial_s^{i_2}k \star \ldots \star \partial_s^{i_n}k$ where $m=i_1 + \ldots+i_n$ is the total number of derivatives.
 
 The following interpolation inequality for closed curves appears in
 \cite{DKS02}.
 
 \begin{interp} \label{T:interp}
 Let $\gamma: I \rightarrow \mathbb{R}^2$ be a smooth closed curve.  Then for any term $P_n^m\left( k\right)$ with $n\geq 2$ that contains derivatives of $k$ of order at most $\ell-1$,
$$\int_I \left| P_n^m\left( k \right)\right| ds \leq c \, L^{1-m-n} \left\| k \right\|_2^{n-p} \left\| k \right\|_{\ell, 2}^{p}$$  
where $p = \frac{1}{\ell} \left( m + \frac{1}{2} n - 1\right)$ and $c=c\left( \ell, m, n\right)$.  Moreover, if $m+ \frac{n}{2} < 2\ell+1$ then $p<2$ and for any $\varepsilon>0$,
\begin{equation*}
  \int_I  \left| P_n^m\left( k \right)\right| ds \leq \varepsilon \int_I \left| \partial_{s^\ell} k \right|^2 ds
   + c\, \varepsilon^{\frac{-p}{2-p}} \left( \int_I \left| k\right|^2 ds\right)^{\frac{n-p}{2-p}} + c\left( \int_I \left| k\right|^2 ds \right)^{m+n-1}\mbox{.}
\end{equation*}
\end{interp}

Note that in the above, $\left\| \cdot \right\|_2$ and $\left\| \cdot \right\|_{m, 2}$ denote scale-invariant norms, for example
$$\left\| k \right\|_2 = \left\| k \right\|_{0, 2} = L^{\frac{1}{2}} \left( \int k^2 ds\right)^{\frac{1}{2}}$$
and
$$\left\| k \right\|_{1, 2} = L^{\frac{1}{2}} \left( \int k^2 ds\right)^{\frac{1}{2}} + L^{\frac{3}{2}} \left( \int k_s^2 ds\right)^{\frac{1}{2}} \mbox{;}$$
of course under our particular flow the $L=L_0$ is constant.  With the exception of the statement of this lemma and above discussion, in this paper we will use the notation $\left\| \cdot \right\|_2$ to denote the regular unscaled norms, pointing out explicit scaling factors where relevant.  In our estimates we will also allow the constants $c$ to vary from line to line where they depend only on absolute quantities like $n$, $m$ and, for this flow $L_0$.  Of course when $L_0$ is embedded into our constants it is no longer possible to track scaling through the estimates.\\

Further, we will need the following elementary inequality that can be used to establish a suitable bound on $h$ in terms of the $L^2$ norms of $k$ and of $k_{s^n}$ for each $n\in \mathbb{N}$.  Here and throughout $k_{s^n}^2$ means $\left( k_{s^n}\right)^2$ etc, where $k_{s^n}$ is the $n$-th iterated derivative of $k$ with respect to arclength.

\begin{interp2} \label{T:interp2}
For $n\in \mathbb{N}$, 
$$\int k_{s^{n-1}}^2 ds \leq \left( \int k^2 ds \right)^{\frac{1}{n}} \left( \int k_{s^n}^2 ds \right)^{\frac{n-1}{n}} \mbox{.}$$
\end{interp2}

\noindent \textbf{Proof:} We proceed by induction.  The statement is trivial for $n=1$.  So assume that 
\begin{equation} \label{E:Ind}
  \int k_{s^{i-1}}^2 ds \leq \left( \int k^2 ds \right)^{\frac{1}{i}} \left( \int k_{s^i}^2 ds \right)^{\frac{i-1}{i}}
\end{equation}
and use this to show
$$\int k_{s^{i}}^2 ds \leq \left( \int k^2 ds \right)^{\frac{1}{i+1}} \left( \int k_{s^{i+1}}^2 ds \right)^{\frac{i}{i+1}} \mbox{.}$$

By integration by parts and the H\"{o}lder inequality we have
$$\int k_{s^{i}}^2 ds =-\int k_{s^{i-1}} k_{s^{i+1}} ds \leq \left( \int k_{s^{i-1}}^2 ds \right)^{\frac{1}{2}}  \left( \int k_{s^{i+1}}^2 ds \right)^{\frac{1}{2}}\mbox{.}$$
Inserting on the right hand side the inductive hypothesis \eqref{E:Ind} we obtain
$$\int k_{s^{i}}^2 ds \leq \left( \int k^2 ds\right)^{\frac{1}{2i}} \left( \int k_{s^i}^2 ds\right)^{\frac{i-1}{2i}}  \left( \int k_{s^{i+1}}^2 ds \right)^{\frac{1}{2}}$$
in other words
$$\left( \int k_{s^{i}}^2 ds \right)^{\frac{i+1}{2i}} \leq \left( \int k^2 ds\right)^{\frac{1}{2i}}  \left( \int k_{s^{i+1}}^2 ds \right)^{\frac{1}{2}}$$
which implies
$$\int k_{s^{i}}^2 ds \leq \left( \int k^2 ds \right)^{\frac{1}{i+1}} \left( \int k_{s^{i+1}}^2 ds \right)^{\frac{i}{i+1}}$$
as required.\hspace*{\fill}$\Box$

\begin{hbound} \label{T:hbound}
For each $n\in \mathbb{N}$, the global term $h\left( t\right)$ may be estimated as 
$$\left| h\left( t\right) \right| \leq \frac{1}{2\pi} \left( \int k^2 ds\right)^{1-\frac{1}{n}} \left( \int k_{s^n}^2 ds\right)^{\frac{1}{n}} \mbox{.}$$
\end{hbound}

\noindent \textbf{Proof:} Again proceed by induction.  The statement is trivially true for $n=1$ from the definition of $h\left( t\right)$.  Assume then that
\begin{equation} \label{E:ind2}
  h\left( t\right) \leq \frac{1}{2\pi} \left( \int k^2 ds\right)^{1-\frac{1}{i}} \left( \int k_{s^i}^2 ds\right)^{\frac{1}{i}}
\end{equation}
and use this to show
$$h\left( t\right) \leq \frac{1}{2\pi} \left( \int k^2 ds\right)^{1-\frac{1}{i+1}} \left( \int k_{s^{i+1}}^2 ds\right)^{\frac{1}{i+1}} \mbox{.}$$
From Lemma \ref{T:interp2} we have
$$\int k_{s^{i}}^2 ds \leq \left( \int k^2 ds \right)^{\frac{1}{i+1}} \left( \int k_{s^{i+1}}^2 ds \right)^{\frac{i}{i+1}} \mbox{.}$$
Substituting this into \eqref{E:ind2} we obtain
$$h\left( t\right) \leq \frac{1}{2\pi} \left( \int k^2 ds\right)^{1-\frac{1}{i}} \left[ \left( \int k^2 ds \right)^{\frac{1}{i+1}} \left( \int k_{s^i}^2 ds \right)^{\frac{i}{i+1}} \right]^{\frac{1}{i}}$$
which simplifies to the required expression.\hspace*{\fill}$\Box$
\mbox{}\\

The following evolution equations for various geometric quantities under the flow \eqref{E:theflow} will be used in our analysis.  These are easily derived similarly as in \cite{W13}, for example.

\begin{evlneqns} \label{T:evlneqns}
Under the flow \eqref{E:theflow}, 
\begin{enumerate}
  \item[\textnormal{(i)}] $\frac{d}{dt} L\left[ \gamma\right] = 0$;
  \item[\textnormal{(ii)}] $\frac{d}{dt} A\left[ \gamma\right] = -h\left( t\right) L_0$;
  \item[\textnormal{(iii)}] $\frac{d}{dt} \int k^2 ds =-2 \int k_{ss}^2 ds + 3 \int k^2 k_s^2 ds + h\left( t\right) \int k^3 ds$;
  \item[\textnormal{(iv)}] 
\begin{multline*}
  \frac{d}{dt} K_{\mbox{osc}} = -2 L_0 \int k_{ss}^2 ds + 3 L_0 \int \left( k-\overline{k}\right)^2 k_s^2 ds + 6 L_0 \overline{k} \int \left( k- \overline{k} \right) k_s^2 ds  \\
  + 2 \overline{k}^2 L_0 \int k_s^2 ds + L_0 h\left( t\right) \left[ \int \left( k - \overline{k}\right)^3 ds + 3 \overline{k} \int \left( k - \overline{k}\right)^2 ds \right]\mbox{;}
  \end{multline*}
  \item[\textnormal{(v)}]
  $\frac{d}{dt} \int k_s^2 ds = -2 \int k_{s^3}^2 ds + 2 \int k^2 k_{ss}^2 ds + \frac{1}{3} \int k_s^4 ds + 5 h\left( t\right) \int k \, k_s^2 ds$;
  \item[\textnormal{(vi)}]
  $\frac{d}{dt} \int k_{ss}^2 ds = -2 \int k_{s^4}^2 ds + 2 \int k^2 k_{s^3}^2 ds - \int k_s^2 k_{ss}^2 ds + 7 h\left( t\right) \int k\, k_{ss}^2 ds$;
  \end{enumerate}
  Moreover, for $m\in \mathbb{N}\cup \left\{ 0 \right\}$,
  \begin{enumerate}
  \item[\textnormal{(vii)}]
  $\frac{d}{dt} \int k_{s^m}^2 ds = - 2 \int k_{s^{m+2}}^2 ds + \int k_{s^m} P_3^{m+2}\left( k \right) ds + h\left( t\right) \int k_{s^m} P_2^m\left( k\right) ds$.
\end{enumerate}
\end{evlneqns}

Here $L=L_0$ is the constant length of the evolving curve $\gamma_t$.\\

We complete this section with a statement of short-time existence for solutions of \eqref{E:theflow}.  

\begin{ste}
Suppose $\gamma_0:\mathbb{S}^1 \rightarrow \mathbb{R}^2$ is a regular curve parametrised by arclength of class $C^1 \cap W^{2,2}$ with $\left\| k \right\|_2 < \infty$.  Then there exists a maximal $T\in \left( 0, \infty \right]$ and a one-parameter family of immersions $\gamma: \mathbb{S}^1 \times \left[ 0, T\right) \rightarrow \mathbb{R}^2$ parametrised by arclength, smooth for $t>0$ and solving \eqref{E:theflow} with $\gamma\left( \cdot, 0 \right) = \gamma_0$.  The family is unique up to the group of invariances of the equation \eqref{E:theflow}.
\end{ste}

\section{The oscillation of curvature} \label{S:estimates}

As a first step, using the isoperimetric inequality, Lemma \ref{T:PSW} and Lemma \ref{T:evlneqns} we can see that $K_{\mbox{osc}}$ is an $L^1$ function in time:

\newtheorem{KoscL1}{Lemma}[section]
\newtheorem{Koscbound}[KoscL1]{Lemma}
\newtheorem{Koscbound2}[KoscL1]{Proposition}

\begin{KoscL1} \label{T:KoscL1}
If $\gamma : \mathbb{S}^1\times \left[ 0, T\right) \rightarrow \mathbb{R}^2$ solves \eqref{E:theflow}, then
$$\left\| K_{\mbox{osc}} \right\|_1 \leq \frac{L_0^2}{2\pi} \left( \frac{L_0}{4\pi} - A_0 \right)$$
where $A_0$ denotes the signed enclosed area of $\gamma\left( \cdot, 0\right)$.
\end{KoscL1}

\noindent \textbf{Proof:} We estimate
$$K_{\mbox{osc}} = L_0 \int \left( k - \overline{k}\right)^2 ds \leq \frac{L_0^3}{4\pi^2} \left\| k_s \right\|_2^2 = \frac{L_0^2}{2\pi} \frac{dA}{dt} \mbox{.}$$
It follows that
$$\left\| K_{\mbox{osc}} \right\|_1 \leq \frac{L_0^2}{2\pi} \left( A\left( t\right) - A_0 \right) \leq \frac{L_0^2}{2\pi} \left( \frac{L_0^2}{4\pi} - A_0 \right)$$
where we have written $A\left( t\right) = A\left[ \gamma_t\right]$.\hspace*{\fill}$\Box$
\mbox{}\\[8pt]

With a little extra work we may obtain a pointwise estimate while $K_{\mbox{osc}}$ is small.  Compared with \cite{W13} we lose some terms since $L_0$ is constant but we pick up two others from $h\left( t\right)$, one of which needs to be estimated.

\begin{Koscbound} \label{T:Koscbound}
Suppose $\gamma : \mathbb{S}^1\times \left[ 0, T\right) \rightarrow \mathbb{R}^2$ solves \eqref{E:theflow}.  If there exists a $T^*$ such that for $t\in \left[ 0, T^*\right)$ we have
$$K_{\mbox{osc}}\left( t\right) \leq 2 K^*$$
then during this time we have
$$K_{\mbox{osc}}\left( t\right) \leq K_{\mbox{osc}}\left( 0\right) + \frac{16\pi^3 \omega^3}{L_0^2} \left( A\left( t\right) - A_0 \right) \mbox{.}$$
\end{Koscbound}

\noindent \textbf{Proof:} Beginning with the evolution equation of Lemma \ref{T:evlneqns} (iv), we estimate using Lemma \ref{T:PSW} as in \cite{W13}
$$3 L_0 \int \left( k - \overline{k}\right)^2 k_s^2 ds \leq \frac{3 L_0}{2\pi} K_{\mbox{osc}} \left\| k_{ss} \right\|_2^2$$
and 
$$6 L_0 \overline{k} \int \left( k - \overline{k} \right) k_s^2 ds \leq 6 \, \omega L_0 K_{\mbox{osc}}^{\frac{1}{2}} \left\| k_{ss}\right\|_2^2 \mbox{.}$$
While we may neglect the $h\left( t\right) \int \left( k- \overline{k} \right)^2 ds$ term, we estimate the other $h\left( t\right)$ term as follows:
$$\int k_s^2 ds = - \int k\, k_{ss} ds = -\int \left( k - \overline{k}\right) k_{ss} ds \leq \left( \int \left( k-\overline{k}\right)^2 ds \right)^{\frac{1}{2}} \left( \int k_{ss}^2 ds \right)^{\frac{1}{2}}$$
and
$$\int \left( k - \overline{k}\right)^3 ds \leq \left\| k - \overline{k} \right\|_\infty \int \left( k - \overline{k}\right)^2 ds$$
thus
$$L_0 \, h\left( t\right) \int \left( k - \overline{k}\right)^3 ds \leq \frac{L_0}{2\pi \omega} \left[  \int \left( k - \overline{k}\right)^3 ds \right]^{\frac{3}{2}} \left\| k - \overline{k} \right\|_\infty \left\| k_{ss} \right\|_2^2 \leq \frac{L_0}{4\pi^2 \sqrt{2\pi} \omega} K_{\mbox{osc}}^{\frac{3}{2}} \left\| k_{ss} \right\|_2^2 \mbox{.}$$
Using also Lemma \ref{T:evlneqns} (ii) we obtain
\begin{equation} \label{E:Kosc}
  \frac{d}{dt} K_{\mbox{osc}} + L_0 \left( 2 - \frac{1}{4\pi^2 \sqrt{2\pi} \, \omega} K_{\mbox{osc}}^{\frac{3}{2}} - \frac{3}{2\pi} K_{\mbox{osc}} - 6 \omega K_{\mbox{osc}}^{\frac{1}{2}} \right) \int k_{ss}^2 ds \leq \frac{16\pi^3 \omega^3}{L_0^2} \frac{dA}{dt} \mbox{.}
\end{equation}
Provided the coefficient of $\int k_{ss}^2 ds$ remains positive on the interval $\left[ 0, T^*\right)$, that is, we take $2K^*$ as the smallest positive solution of 
$$2 - \frac{1}{4\pi^2 \sqrt{2\pi} \, \omega} K_{\mbox{osc}}^{\frac{3}{2}} - \frac{3}{2\pi} K_{\mbox{osc}} - 6 \, \omega K_{\mbox{osc}}^{\frac{1}{2}}=0$$
then we obtain the result by integration in time.\hspace*{\fill}$\Box$
\mbox{}\\

Working similarly as in \cite{W13} in the case $\omega=1$ we next show that if initially $K_{\mbox{osc}}$ is sufficiently small and the isoperimetric ratio $I\left( 0\right)= \frac{L_0^2}{4\pi A\left( 0\right)}$ is sufficiently close to $1$ then $K_{\mbox{osc}}$ remains under control under \eqref{E:theflow} as long as the solution exists.  When $\omega=1$, we may estimate from the expression in the previous proof that $2K^* \approx 0.09$.

\begin{Koscbound2} \label{T:Koscbound2}
Suppose $\gamma : \mathbb{S}^1\times \left[ 0, T\right) \rightarrow \mathbb{R}^2$ has $\omega=1$ and solves \eqref{E:theflow}.  Then if
$$K_{\mbox{osc}}\left( 0 \right) < K^* \mbox{ and } I\left( 0\right) < \frac{4\pi^2}{4\pi^2 - K^*}$$
we have that $K_{\mbox{osc}} < 2K^*$ for all $t\in \left[ 0, T\right)$.
\end{Koscbound2}

\noindent \textbf{Proof:} Suppose for the sake of establishing a contradiction there exists a first time $T^* <T$ for which $K_{\mbox{osc}}\left( t\right) = 2 K^*$.  In view of Lemma \ref{T:Koscbound} we have
$$K_{\mbox{osc}}\left( T^* \right) \leq K_{\mbox{osc}}\left( 0\right) + \frac{16 \pi^3}{L_0^2} \left( A\left( T^*\right) - A\left( 0 \right) \right) \leq K_{\mbox{osc}}\left( 0 \right) + 4 \pi^2 \left( 1 - \frac{1}{I\left( 0 \right)} \right) < K^* + K^* = 2 K^*$$
a contradiction.  We conclude that $K_{\mbox{osc}} < 2K^*$ for all $t\in \left[ 0, T\right)$.\hspace*{\fill}$\Box$

\section{Global existence} \label{S:global}
\newtheorem{blowup}{Theorem}[section]
\newtheorem{longtime}[blowup]{Corollary}
\newtheorem{ks}[blowup]{Proposition}
\newtheorem{decay}[blowup]{Corollary}

Similarly as in \cite{DKS02}[Theorem 3.1] for the curve diffusion flow, we first show that if the maximal existence time is finite, then the curvature must blow up in $L^2$.

\begin{blowup} \label{T:blowup}
Let $\gamma : \mathbb{S}^1\times \left[ 0, T\right) \rightarrow \mathbb{R}^2$ be a maximal solution of \eqref{E:theflow}.  If $T< \infty$ then
$$\int k^2 ds \geq c\left( T - t\right)^{-\frac{1}{4}} \mbox{.}$$
\end{blowup}

\noindent \textbf{Proof:} We proceed similarly as in \cite{DKS02} and make the necessary adjustments.   Integrating once by parts on the second term on the right hand side of Lemma \ref{T:evlneqns}, (vii) we can ensure the highest derivative appearing is $k_{s^{m+1}}$ and thus using Lemma \ref{T:interp} we have as in \cite{DKS02}
$$\int k_{s^m} P_3^{m+2}\left( k\right) ds \leq \varepsilon \int k_{s^{m+2}}^2 ds + c_m\left( \varepsilon\right) \left( \int k^2 ds\right)^{2m+5} \mbox{.}$$

For the $h\left( t\right)$ term we work as follows.  Using Lemma \ref{T:interp} we have
\begin{multline*}
  \int k_{s^m} P_2^m\left( k\right) ds \leq c\left( L_0\right) \left( \int k^2 ds\right)^{\frac{m+ \frac{5}{2}}{2m+2}}\left[ \left( \int k^2 ds\right)^{\frac{2m+ \frac{1}{2}}{2m+2}} + \left( \int k_{s^{m+1}}^2 ds \right)^{\frac{2m+\frac{1}{2}}{2m+2}} \right]\\
  = c\left( L_0 \right) \left[ \left( \int k^2 ds\right)^{\frac{3}{2}} + \left( \int k^2 ds\right)^{\frac{m+ \frac{5}{2}}{2m+2}}\left( \int k_{s^{m+1}}^2 ds \right)^{\frac{2m+\frac{1}{2}}{2m+2}} \right]
  \end{multline*}
where we have used Lemma \ref{T:PSW} to estimate all the intermediate derivatives in $\left\| k \right\|_{m+1, 2}$ by $\int k_{s^{m+1}}^2 ds$.  Applying now Lemma \ref{T:interp2} on the last term we obtain
$$ \int k_{s^m} P_2^m\left( k\right) ds \leq c\left( L_0 \right) \left[ \left( \int k^2 ds\right)^{\frac{3}{2}} + \left( \int k^2 ds\right)^{\frac{2m+11}{4m+8}} \left( \int k_{s^{m+2}}^2 ds \right)^{\frac{2m+\frac{1}{2}}{2m+4}} \right] \mbox{.}$$
Using also from Lemma \ref{T:hbound}
$$h\left( t\right) \leq \frac{1}{2\pi} \left( \int k^2 ds\right)^{\frac{m+1}{m+2}} \left( \int k_{s^{m+2}}^2 ds\right)^{\frac{1}{m+2}}$$
we have
\begin{multline*}
  h\left( t\right) \int k_{s^m} P_2^m\left( k\right) ds \\
  \leq c\left( L_0\right) \left[ \left( \int k^2 ds\right)^{\frac{5m+8}{2m+4}} \left( \int k_{s^{m+2}}^2 ds \right)^{\frac{1}{m+2}} + \left( \int k^2 ds\right)^{\frac{6m+15}{4m+8}} \left( \int k_{s^{m+2}}^2 ds \right)^{\frac{2m + \frac{5}{2}}{2m+4}} \right]\mbox{.}
\end{multline*}
Importantly, above the powers of $\int k_{s^{m+2}}^2 ds$ are less than $1$.  We now estimate each of the terms on the right hand side using Young's inequality to obtain
$$ h\left( t\right) \int k_{s^m} P_2^m\left( k\right) ds \leq \varepsilon \int k_{s^{m+2}}^2 ds + c \left( \int k^2 ds \right)^{\frac{5m+8}{2m+2}} + c \left( \int k^2 ds\right)^{2m+5} \mbox{.}$$
Substituting these estimates into Lemma \ref{T:evlneqns}, (vii) we obtain for suitably small $\varepsilon$
\begin{equation} \label{E:dt}
  \frac{d}{dt} \int k_{s^m}^2 ds + \int k_{s^{m+2}}^2 ds \leq c \left( \int k^2 ds \right)^{\frac{5m+8}{2m+2}} + c \left( \int k^2 ds\right)^{2m+5} \mbox{.}
  \end{equation}

The proof may now be completed similarly as in \cite{DKS02} restricting to the codimension $1$ case (see also \cite{W15}, for example).  The idea is that if on the contrary $T< \infty$ is maximal but $\int k^2 ds \leq \Lambda <\infty$ for all $t<T$, then the flow can be smoothly extended beyond $t=T$ via short-time existence, contradicting the maximality of $T$.  Hence it must be the case that $\limsup_{t\rightarrow T} \int k^2 ds \rightarrow \infty$ if $T$ is finite.

To determine the blow up rate observe that for $m=0$ we have
$$\frac{d}{dt} \int k^2 ds \leq c  \left( \int k^2 ds\right)^4 + c \left( \varepsilon \right) \left( \int k^2 ds\right)^5 \mbox{.}$$
Since $\int k^2 ds$ blows up as $t\rightarrow T$, the power $5$ term dominates leading to the given blow-up rate by solving the ordinary differential inequality. \hspace*{\fill}$\Box$\\

\noindent \textbf{Remark:} The inequality \eqref{E:dt} above implies under the condition $\int k^2 ds \leq \Lambda$ that all curvature derivatives are bounded in $L^2$.  This fact is used later together with integration by parts and the H\"{o}lder inequality to obtain exponential convergence of higher curvature derivatives.\\

\begin{longtime} \label{T:longtime}
Suppose $\gamma : \mathbb{S}^1\times \left[ 0, T\right) \rightarrow \mathbb{R}^2$ solves \eqref{E:theflow} and satisfies the conditions of Theorem \ref{T:main}.  Then $T=\infty$.
\end{longtime}

\noindent \textbf{Proof:} Suppose on the contrary that $\gamma$ satisfies the conditions of Proposition \ref{T:Koscbound2} and $T<\infty$.  We know from Theorem \ref{T:blowup} that $\left\| k\right\|_2^2 \rightarrow \infty$ as $t\rightarrow T$.  Then
$$K_{\mbox{osc}} = L_0 \int \left( k - \overline{k}\right)^2 ds = L_0 \left\| k \right\|_2^2 - 2 \pi^2 \rightarrow \infty$$
as $t\rightarrow T$.  But this contradicts Lemma \ref{T:Koscbound2}.  We conclude that it must be the case that $T=\infty$.\hspace*{\fill}$\Box$\\

It remains to classify the limit as $t\rightarrow \infty$ where $\gamma$ satisfies the conditions of Proposition \ref{T:Koscbound2}.  By Lemma \ref{T:KoscL1}, $K_{\mbox{osc}} \in L^1\left( \left[ 0, \infty \right) \right)$ so if we show $K_{\mbox{osc}}'$ is bounded we may conclude $K_{\mbox{osc}} \rightarrow 0$ and therefore limit curves are circles of circumference length $L_0$.  In view of \eqref{E:Kosc}, Lemma \ref{T:evlneqns}, (ii) and \eqref{E:h}, to show $K_{\mbox{osc}}'$ is bounded amounts to establishing a uniform bound for $\left\| k_s\right\|_2^2$. 

\begin{ks} \label{T:ks}
Suppose $\gamma : \mathbb{S}^1\times \left[ 0, \infty\right) \rightarrow \mathbb{R}^2$ solves \eqref{E:theflow} and satisfies the conditions of Theorem \ref{T:main}.  Then there exists a constant $c_1>0$ depending only on $\gamma_0$ such that
$$\left\| k_s\right\|_2^2 \leq c_1 \mbox{.}$$
\end{ks}

\noindent \textbf{Proof:} We use Lemma \ref{T:evlneqns} (v) and estimate terms in terms of $K_{\mbox{osc}}$ whose behaviour is under control.  Using integration by parts
$$\int k_s^4 ds =-3\int k\, k_s^2 k_{ss} ds \leq \frac{1}{2} \int k_s^4 ds + \frac{9}{2} \int k^2 k_{ss}^2 ds$$
so
$$\int k_s^4 ds \leq 9 \int k^2 k_{ss}^2 ds$$
and the positive non-$h\left( t\right)$ terms may be estimated as in \cite{W13}, page 944:
$$\int k^2 k_{ss}^2 ds \leq \left( \frac{1}{27} + \frac{K_{\mbox{osc}}}{\pi}\right) \left\| k_{s^3}\right\|_2^2+ \frac{27^3 \overline{k}^8 L_0 K_{\mbox{osc}}}{2\pi^2} \mbox{.}$$
For the $h\left( t\right)$ term we estimate using integration by parts
$$\int k_s^2 ds = - \int \left( k - \overline{k} \right) k_{ss} ds \leq \left[ \int \left( k - \overline{k} \right)^2 ds \right]^{\frac{1}{2}} \left( \int k_{ss}^2 ds \right)^{\frac{1}{2}} = \frac{K_{\mbox{osc}}^{\frac{1}{2}}}{L_0^{\frac{1}{2}}} \left\| k_{ss}\right\|_2$$
and
\begin{align*}
  \int k\, k_s^2 ds &= - \frac{1}{2} \int k^2 k_{ss} ds = - \frac{1}{2} \int \left( k - \overline{k}\right)^2 k_{ss} ds - \overline{k} \int \left( k - \overline{k} \right) k_{ss} ds \\
  &\leq \frac{K_{\mbox{osc}}}{2 L_0} \left\| k_{ss} \right\|_{\infty} + \overline{k} \left[ \int \left( k - \overline{k} \right)^2 ds \right]^{\frac{1}{2}} \left[ \int k_{ss}^2 ds \right]^{\frac{1}{2}}\\
  &\leq \frac{1}{2 L_0^{\frac{1}{2}}} \left( \frac{K_{\mbox{osc}}}{\sqrt{2\pi}} + 2 K_{\mbox{osc}}^{\frac{1}{2}}\right) \left\| k_{s^3} \right\|_2 \mbox{.}
  \end{align*}
  It follows that
  $$5\, h\left( t\right) \int k \, k_s^2 ds \leq \frac{5 \, K_{\mbox{osc}}}{8 \pi^2} \left( \frac{K_{\mbox{osc}}^{\frac{1}{2}}}{\sqrt{2\pi}} + 2 \right) \left\| k_{s^3}\right\|_2^2 \mbox{.}$$
Substituting these estimates into Lemma \ref{T:evlneqns} (v) we obtain
\begin{equation*}
  \frac{d}{dt} \int k_s^2 ds + \left[ 2 - 5\left( \frac{1}{27} + \frac{K_{\mbox{osc}}}{\pi} \right) - \frac{5 \, K_{\mbox{osc}}}{8 \pi^2} \left( \frac{K_{\mbox{osc}}^{\frac{1}{2}}}{\sqrt{2\pi}} + 2 \right) \right] \left\| k_{s^3}\right\|_2^2
  \leq \frac{5 \times 27^3 \overline{k}^8 L_0 K_{\mbox{osc}}}{2\pi^2} \mbox{.}
  \end{equation*}
The smallness of $K_{\mbox{osc}}$ ensures the coefficient of $\left\| k_{s^3}\right\|_2^2$ is positive (about $1.642$) and the right hand side is bounded.  The proof is then completed using Lemma \ref{T:PSW} and integrating.\hspace*{\fill}$\Box$ 
\mbox{}\\

We are now able to estimate length of the set of times for which $k$ is not strictly positive.\\

\noindent \textbf{Proof of Theorem \ref{T:pos}}: Following the idea in \cite{W13} to rearrange $\gamma$ in time if necessary, we may assume that $k\left( \cdot, t\right) \not > 0$ for all $t\in \left[ 0, t_0\right)$ while $k\left( \cdot, t\right) >0$ for all $t\in \left[ t_0, \infty \right)$.  Let us suppose for the sake of establishing a contradiction that $t_0 >\frac{L_0^2}{4\, \pi^3}\left( \frac{L_0^2}{4\, \pi} - A_0 \right)$.

Working in the interval $\left[ 0, t_0 \right)$ we know that $k$ has a zero therefore
$$\int k^2 ds \leq \frac{L_0^2}{2\pi^2} \int k_s^2 ds \mbox{.}$$
Inserting this into Lemma \ref{T:evlneqns} (ii) we estimate for almost every $t\in \left[ 0, t_0\right)$
$$\frac{d}{dt} A = \frac{L_0}{2\pi}\int k_s^2 ds \geq \frac{\pi}{L_0} \int k^2 ds \geq \frac{4\pi^3}{L_0^2}$$
where in the last step we used the H\"{o}lder inequality.  It follows that for almost every $t\in \left[ 0, t_0\right)$,
$$A\left( t\right) \geq A_0 + \frac{4\pi^3}{L_0^2}t \mbox{.}$$
Taking the limit $t\rightarrow t_0$ we establish a contradiction to the isoperimetric inequality in view of our assumption on $t_0$.\hspace*{\fill}$\Box$
\mbox{}\\

To complete the proof of Theorem \ref{T:main} it remains to show that the limit circle is unique and convergence is exponential.  Since we have convergence to circles of radius $\frac{L_0}{2\pi}$,  we can be sure that $k\left( x, t\right) \in \left[ \frac{\pi}{L_0}, \frac{3\pi}{L_0} \right]$ say for all $t\geq t_1$.  For such times we also have $\left\| k \right\|_\infty \leq \frac{3\pi}{L_0}$.

\begin{decay} \label{T:decay}
Suppose $\gamma: \mathbb{S}^1 \times \left[ 0, \infty \right) \rightarrow \mathbb{R}^2$ solves \eqref{E:theflow} and satisfies the assumptions of Theorem \ref{T:main}.  Then there exists constants $c_1, c_2 >0$ depending only on $\gamma_0$ such that
$$\left\| k_{ss} \right\|_2^2 \leq c_1\, \mbox{e}^{-c_2 t} \mbox{.}$$
\end{decay}

\noindent \textbf{Proof:} For all $t \geq t_1$, the curve $\gamma\left( \cdot, t\right)$ is convex so the $h\left( t\right)$ term of Lemma \ref{T:evlneqns}, (vi), is negative.  For the nonnegative term we estimate
\begin{multline*}
  \int k^2 k_{s^3}^2 ds = \int k_{s^3}^2 \left( k - \overline{k} \right)^2 ds + 2\, \overline{k} \int k\, k_{s^3}^2 ds - \overline{k}^2 \int k_{s^3}^2 ds\\
    \leq \int k_{s^3}^2 \left( k - \overline{k} \right)^2 ds + \frac{1}{2} \int k^2 k_{s^3}^2 ds + \overline{k}^2 \int k_{s^3}^2 ds \mbox{;}
   \end{multline*}
   absorbing on the left yields
\begin{equation} \label{E:k2k32}
  \int k^2 k_{s^3}^2 ds \leq 2 \int k_{s^3}^2 \left( k- \overline{k}\right)^2 ds + 2\, \overline{k}^2 \int k_{s^3}^2 ds \mbox{.}
  \end{equation}
   Now for $\sigma_1, \sigma_2 >0$,
   $$\int k_{s^3}^2 ds = - \int k_{ss} k_{s^4} ds \leq \sigma_1 \int k_{s^4}^2 ds + \frac{1}{4\sigma_1} \int k_{ss}^2 ds \mbox{,}$$
   \begin{multline*}
     \int k_{ss}^2 ds = - \int k_{s^3} k_s ds = \int k_{s^4} k ds = \int k_{s^4} \left( k- \overline{k}\right) ds\\
     \leq \sigma_2 \int k_{s^4}^2 ds + \frac{1}{4\sigma_2} \int \left( k - \overline{k}\right)^2 ds = \sigma_2 \int k_{s^4} ds + \frac{1}{4\sigma_2} \frac{K_{\mbox{osc}}}{L_0}
\end{multline*}
and using Lemma \ref{T:PSW} (ii),
$$\int k_{s^3}^2 \left( k - \overline{k}\right)^2 ds \leq \left\| k_{s^3}\right\|_\infty^2 \frac{K_{\mbox{osc}}}{L_0} \leq \frac{K_{\mbox{osc}}}{2\pi} \int k_{s^4}^2 ds \mbox{.}$$
Substituting into \eqref{E:k2k32} yields
$$\int k^2 k_{s^3}^2 ds \leq \left[ 2 \overline{k}^2 \left( \sigma_1 + \frac{\sigma_2}{4\sigma_1} \right) + \frac{K_{\mbox{osc}}}{\pi} \right] \int k_{s^4}^2 ds + \frac{\overline{k}^2}{8 \sigma_1 \sigma_2} \frac{K_{\mbox{osc}}}{L_0}$$
and thus from Lemma \ref{T:evlneqns} (vi) we have
$$\frac{d}{dt} \int k_{ss}^2 ds + 2 \int k_{s^4}^2 ds \leq \left[ 4 \overline{k}^2 \left( \sigma_1 + \frac{\sigma_2}{4\sigma_1} \right) + \frac{2\, K_{\mbox{osc}}}{\pi} \right] \int k_{s^4}^2 ds + \frac{\overline{k}^2}{4 \sigma_1 \sigma_2} \frac{K_{\mbox{osc}}}{L_0} \mbox{.}$$
Since $K_{osc} \rightarrow 0$ we can choose $\sigma_1$ and $\sigma_2$ small enough to obtain for a small $\delta>0$ and all $t\geq t_2 \geq t_1$ that
$$\frac{d}{dt} \int k_{ss}^2 ds \leq - \delta \int k_{s^4}^2 ds + \frac{\overline{k}^2}{4 \sigma_1 \sigma_2} \frac{K_{\mbox{osc}}}{L_0} \mbox{.}$$
Applying Lemma \ref{T:PSW} (i) twice and recalling Lemma \ref{T:KoscL1} we obtain the desired exponential convergence all $t\geq t_2$ where $t_2\geq t_1$ using Gr\"{o}nwall's inequality.  The result can be extended to $\left[ 0, \infty \right)$ by adjusting the constants.\hspace*{\fill}$\Box$ 
\mbox{}\\

\noindent \textbf{Completion of the proof of Theorem \ref{T:main}}:  Using Corollary \ref{T:decay} with Lemma \ref{T:PSW} gives in turn exponential decay of $\int k_s^2 ds$ (hence decay of $h\left( t\right)$) and $\int \left( k - \overline{k}\right)^2 ds$ and the corresponding $L^{\infty}$ norms.  This implies subconvergence of the flow to circles with perimeter length $L_0$.  Exponential decay of the higher curvature derivatives follows by interpolation using the uniform bounds on $\int k_{s^m}^2 ds$ in the proof of Theorem \ref{T:blowup} that apply since $\int k^2 ds$ is uniformly bounded.  A stability argument may be used to obtain stronger convergence to the circle.  In view of uniform convergence of $k$, we may work from a time beyond which $\gamma\left( \cdot, t\right)$ remains convex.  Then we may use the radial graph parametrisation $\gamma\left( z, t\right) = \rho\left( z, t\right) z$ for $z\in \mathbb{S}^1$ and consider the length preserving deformations $\rho\left( z, t\right) = \rho_\infty + \varepsilon u\left( z, t\right)$, where $\rho_\infty = \frac{L_0}{2\pi}$ is the radius of the limiting circle.  That length is preserved corresponds to the condition $\int_\gamma u\, dz =0$.  In the radial graph parametrisation the linearised operator is
$$\mathcal{L}u= \rho_\infty^{-4} \left( u_{x^4} + u_{xx} \right) \mbox{.}$$
All eigenvalues of $\mathcal{L}$ are negative with the exception of a single zero eigenvalues corresponding to translations.  Uniqueness of the limit then follows by Hale-Raugel's convergence theorem \cite{HR92}.  In view of the previous results, all solutions satisfying the conditions of Theorem \ref{T:main} thus converge exponentially to a unique limiting circle. \hspace*{\fill}$\Box$ 
\mbox{}\\[8pt]
\noindent \textbf{Remark:} Exponential decay of the speed allows us to bound the region of the plane in which the solution lies relative to $\gamma_0$ via standard arguments. Specifically we may bound the distance travelled by any point on the initial curve $\gamma_0$ as follows
$$\left| \gamma\left( x, t\right) - \gamma\left( x, 0\right) \right| = \left| \int_0^t \frac{\partial \gamma}{\partial t}\left( x, \tau\right) \right| \leq \int_0^t \left| h\left( \tau\right) - \kappa_{ss} \right| d\tau
\leq \frac{C}{\delta} \left( 1- e^{-\delta t}\right) \leq \frac{C}{\delta}\mbox{.}$$

\section{Self-similar solutions} \label{S:selfsim}

\newtheorem{stny}{Lemma}[section]
A self-similar solution to a curvature flow equation such as \eqref{E:theflow} is a solution whose image maintains the same shape as it evolves; it changes in time only by scaling, translation and/or rotation.  In the present setting the length constraint rules out expanding and contracting self-similar solutions, so we focus here on stationary solutions, translators and rotators.

We begin with the following simple observation:

\begin{stny}
The only smooth, closed stationary solutions to \eqref{E:theflow} are multiply-covered circles.
\end{stny}

\noindent \textbf{Proof:} Such solutions satisfy
$$h\left( t\right) - k_{ss} \equiv 0 \mbox{,}$$
in other words
$$\frac{1}{2\, \pi\, \omega}\int  k_s^2 ds - k_{ss} \equiv 0 \mbox{.}$$
Integrating this equation over $\gamma$ we obtain
$$\frac{L_0}{2\, \pi \omega} \int_\gamma k_s^2 ds \equiv 0$$
so such closed curves have $k_s \equiv 0$ and are thus circles, whose length is controlled via the prescribed $L_0$.\hspace*{\fill}$\Box$
\mbox{}\\[8pt]

In fact, such curves also turn out to be the only possible closed rotators under \eqref{E:theflow}:\\

\noindent \textbf{Proof of Proposition \ref{T:rot8}}: Similarly as in \cite[Section 7]{EGBMWW14} but for the flow \eqref{E:theflow}, curves evolving purely by rotation must satisfy
$$h\left( t\right) - k_{ss} = 2 \, S\left( t\right) \left< \gamma, \gamma_s\right>$$
for some function $S$.  Integrating this equation we obtain
$$\frac{L_0}{2\, \pi} \int_\gamma k_s^2 ds =  S\left( t\right) \int_\gamma \frac{d}{ds} \left| \gamma\right|^2 ds = 0 \mbox{,}$$
regardless of $S\left( t\right)$.  Hence again we must have $k_s \equiv 0$ and thus closed curves are circles.\hspace*{\fill}$\Box$
\mbox{}\\[8pt]

Finally we consider the case of closed curves translating under \eqref{E:theflow}.\\

\noindent \textbf{Proof of Proposition \ref{T:transl8}}: Similarly as in \cite[Section 5]{EGBMWW14}, translators must satisfy
\begin{equation} \label{E:transl8}
  h\left( t\right) - k_{ss} = \left< V, \nu \right>
\end{equation}
for some constant vector $V$ in the direction of translation.  Integrating this equation we obtain
$$\frac{L_0}{2\, \pi} \int_\gamma k_s^2 ds = \int_\gamma \left< V, \nu \right> ds = 0$$
regardless of $V$.  Again it follows that $k_s \equiv 0$ and thus closed curves are circles.  It follows now from \eqref{E:transl8} that in fact $V\equiv 0$ so the translators are actually stationary.\hspace*{\fill}$\Box$

\end{document}